\documentclass[10pt]{article}
\usepackage{amssymb,mathrsfs,amsmath,amsthm}

\newcommand{\twostac}[2]{\left(\begin{array}{c}{#1}\\{#2}\end{array}\right)}

\let\a=\alpha
\let\b=\beta

\let\e=\varepsilon

\let\g=\gamma

\let\l=\lambda

\let\m=\mu

\let\N=\nabla
\let\nd=\nabla

\let\r=\rho

\def\bbC{{\mathbb C}}

\def\cL{{\cal L}}

\def\cV{{\cal V}}

\def\gg{{\mathfrak{g}}}

\def\dc{\hbox{$\nabla\!\!\!\! /\;$}}

\let\ptl=\partial

\def\sideremark#1{\ifvmode\leavevmode\fi\vadjust{\vbox to0pt{\vss
 \hbox to 0pt{\hskip\hsize\hskip1em
 \vbox{\hsize2cm\tiny\raggedright\pretolerance10000
 \noindent #1\hfill}\hss}\vbox to8pt{\vfil}\vss}}}%
\def\Lap{\triangle}

\newtheorem{thm}{Theorem}[section]
\newtheorem{rmk}[thm]{Remark}
\newtheorem{cor}[thm]{Corollary}
\newtheorem{lem}[thm]{Lemma}

\title{Conformally Invariant Powers of Spin Operators on the Sphere}
\author{Doojin Hong}

\begin{document}

\maketitle

\begin{abstract}
We give explicit formulas for all odd order differential intertwinors on the subbundle of the bundle of spinor-$k$-forms that are annihilated by the Clifford multiplication over the odd dimensional standard sphere. The Dirac and Rarita-Schwinger operators appear in the case of $k=0$ and $k=1$, respectively.  
\end{abstract}
\section{Introduction}
The Spectrum generating technique \cite{BOO:96} can be applied on the spinor bundle over the sphere $S^n$ with $n$ odd. From the spectral function obtained from the process, one can immediately build conformally invariant powers of the Dirac operator of all odd orders \cite{BO:06}. 

When the same technique is applied on the subbundle of spinor-$k$-forms with $k\ge 1$ that are annihilated by the Clifford multiplication, we get two different spectral functions on two different $K$-type isotypic summands of the section space. The conformally invariant first order differential operator on this subbundle is well known \cite{Branson:99} and it is readily seen that the operator agrees with the spectral function on each isotypic summand. 

However, unlike the Dirac operator case, taking simple powers of this first order operator does not yield higher order conformally invariant operators. We consider a natural operator (generalized gradient) composed with its formal adjoint. This second order operator acts as zero on the one type of isotypic summand and as a nonzero scalar on the other type of isotypic summand. We show that we can build higher order conformally invariant operators from the well-known first order operator together with certain scalar multiples of this second order operator.  That is, we show the formulas for all odd order conformally invariant operators in odd dimensional case. 
\section{Conformally covariant operators}
Let $(M,g)$ be an $n$-dimensional pseudo-Riemannian manifold. 
If $f$ is a (possibly local) diffeomorphism on $M$, we denote by $f\cdot$ the 
natural action of $f$ on tensor fields which acts on vector fields as 
$f\cdot X=(df)X$ and on covariant tensors as $f\cdot \phi=(f^{-1})^*\phi$.  

A vector field $T$ is said to be {\it conformal} with {\it conformal factor} 
$\omega\in C^{\infty}(M)$ 
if 
$$
\mathcal{L}_T\, g=2\omega g\, ,
$$ 
where $\mathcal{L}$ is the Lie derivative. The conformal vector fields form a 
Lie algebra $\mathfrak{c}(M,g)$.
A conformal transformation on $(M,g)$ is a (possibly local) diffeomorphism $h$ 
for which $h\cdot g=\Omega^2 g$ for some positive function $\Omega\in 
C^{\infty}(M)$. The global conformal transformations form a group $\mathscr{C}(M,g)$. 
Let $\mathscr{T}$ be a space of $C^\infty$ tensor fields of some fixed type over $M$. For example, we can take 2-forms or trace-free symmetric covariant three-tensors. 
We have representations \cite{Branson:96} defined by 
\begin{align}\label{int}
\mathfrak{c}(M,g)\stackrel{U_a}{\longrightarrow} \mbox{ End }\mathscr{T},\quad
&U_a(T)=\mathcal{L}_T+a\omega\, \mbox{ and}\\
\mathscr{C}(M,g)\stackrel{u_a}{\longrightarrow} \mbox{ Aut }\mathscr{T}, \quad
&u_a(h)=\Omega^a h\cdot \nonumber
\end{align}
for $a\in \bbC$.

Note that if a conformal vector field $T$ integrates to a one-parameter group of 
global conformal transformation $\{h_{\e}\}$, then
$$
\{U_a(T)\phi\}(x)=\frac{d}{d\e}\Big|_{\e=0}\{u_a(h_{-\e})\phi\}(x)\, .
$$
In this sense, $U_a$ is the infinitesimal representation corresponding to $u_a$. 
\\
A differential operator $D:C^{\infty}(M)\rightarrow C^{\infty}(M)$ is said to be 
{\it infinitesimally conformally covariant of bidegree} $(a,b)$ if 
$$
DU_a(T)\phi=U_b(T)D\phi
$$
for all $T\in \mathfrak{c}(M,g)$ and $D$ is said to be {\it conformally covariant of bidegree} $(a,b)$ if 
$$
Du_a(h)\phi=u_b(h)D\phi
$$
for all $h\in \mathscr{C}(M,g)$. 

To relate conformal covariance to conformal invariance, we recall that 
conformal weight of a bundle $V$ with the induced bundle metric $g_V$ from $g$ is $r$ 
iff 
$$
\bar{g}=\Omega^2 g \Longrightarrow \overline{g_V}=\Omega^{-2r} g_V\, .
$$
Tangent bundle, for instance, has conformal weight -1. Let us denote a bundle $V$ with 
conformal weight $r$ by $V^r$. Then we can impose new conformal weight $s$ on $V^r$ by 
taking tensor product of it with the bundle $I^{(s-r)/n}$ of scalar $((s-r)/n)$-densities \cite{Branson:97-1}. Now if 
we look at an operator of bidegree $(a,b)$ as an operator from the bundle with conformal weight $-a$ to the bundle with conformal weight $-b$, the operator becomes conformally 
invariant.

As an example, let us consider the conformal Laplacian on $M$: 
\begin{equation*}
Y=\triangle+\frac{n-2}{4(n-1)}R, 
\end{equation*}
where $\triangle=-g^{ab}\nabla_a\nabla_b$ and $R$ is the scalar curvature. 
Note that $Y: C^\infty(M) \rightarrow C^\infty(M)$ is conformally covariant of bidegree $((n-2)/2,(n+2)/2)$. That is,
\begin{equation*}
\overline{Y}=\Omega^{-\frac{n+2}{2}}Y\mu(\Omega^{\frac{n-2}{2}})\, ,
\end{equation*}
where $\overline{Y}$ is $Y$ evaluated in $\overline{g}$ and 
$\mu(\Omega^{\frac{n-2}{2}})$ is multiplication by $\Omega^{\frac{n-2}{2}}$.
If we let $V=C^\infty(M)$ and view $Y$ as an operator 
\begin{equation*}
Y: V^{-\frac{n-2}{2}} \rightarrow V^{-\frac{n+2}{2}}\, ,
\end{equation*}
we have, for $\phi\in V^{-\frac{n-2}{2}}$, 
\begin{equation*}
\overline{Y}\,\,\overline{\phi}=\overline{Y\phi}\, ,
\end{equation*}
where $\overline{Y}$, $\overline{\phi}$, and $\overline{Y\phi}$ are $Y$, $\phi$, and $Y\phi$ computed in $\overline{g}$, respectively.
\section{Dominant weights}
Let $\lambda$ be a dominant weight of an irreducible Spin($n$) representation. That is, $\lambda$ is an $l$-tuple $(\l_1,\cdots,\l_l)\in \mathbb{Z}^l\cup(\frac{1}{2}+\mathbb{Z})^l,\, l=[n/2]$, satisfying the inequality constraint (dominant condition)
\begin{equation*}
\begin{array}{ll}
\l_1\ge\cdots\ge\l_l\ge0,&n\text{ odd},\\
\l_1\ge\cdots\ge\l_{l-1}\ge |\l_l|,&n\text{ even}.
\end{array}
\end{equation*}
$\l$ is identified with the highest weight of the irreducible representation of Spin($n$)  \cite{IT:78}. We shall denote by $V(\l)$ the representation with the highest weight $\l$. Those $\lambda\in\mathbb{Z}^l$ are exactly the representations that factor through SO($n$). For example, $V(1,0,\dots,0)$ and $V(1,1,1,0,\dots,0)$ are the defining representation and the three-form representation of SO($n$), respectively and $V(\frac{1}{2},\dots,\frac{1}{2})$ is the the spinor representation in odd dimensional case.

If $M$ is an $n$-dimensional smooth manifold with Spin($n$) structure and $\mathcal{F}$ is the bundle of spin frames, we denote by $\mathbb{V}(\l)$ the associated vector bundle $\mathcal{F}\times_\l V(\l)$. 
\section{Intertwining relation}
Let $G=\text{Spin}_0(n+1,1)$ be the identity component of the $\text{Spin}(n+1,1)$ and $\mathfrak{g}=\mathfrak{k}+\mathfrak{s}$ be a Cartan decomposition of the Lie algebra $\mathfrak{g}$ of $G$. Then, in an Iwasawa decomposition $G=KAN$, the maximal compact subgroup $K$ of $G$ is a copy of $\text{Spin}(n+1)$. Let $M$ be the centralizer of the Lie algebra $\mathfrak{a}$ of $A$ in $K$. Then $M$ is a copy of $\text{Spin}(n)$ and $P=MAN$ is a maximal parabolic subgroup of $G$. Note that $G/P=K/M$ is diffeomorphic to the sphere $S^n$ \cite{BOO:96}. 

Let $V(\lambda)$ be a finite dimensional irreducible representation of $M$. Consider the $G$ module $\mathcal{E}(G;\lambda,\nu)$ of $C^\infty$ functions
\begin{equation*}
F: G\rightarrow V(\lambda) \text{ with }F(gman)=a^{-\nu-\rho}\lambda(m)^{-1}F(g), \, g\in G, m\in M, a\in A, n\in N,
\end{equation*}
where $\rho$ is the half the sum of the positive $(\mathfrak{g},\mathfrak{a})$ roots. This is the space of smooth sections of $\mathbb{V}(\nu,\lambda)$ and is in one-to-one correspondence with the space of smooth sections of $\mathbb{V}(\lambda)$, the $K$ module $\mathcal{E}(K;\lambda)$ of $C^\infty$ functions
\begin{equation*}
f:K\rightarrow V(\lambda) \text{ with }f(km)=\lambda(m)^{-1}f(k),\, k\in K, m\in M.
\end{equation*}
The $K$-finite subspace $\mathcal{E}_K(G;\lambda,\nu)\cong_K \mathcal{E}_K(K;\lambda)$ is defined as 
\begin{equation*}
\bigoplus_{\alpha\in \hat{K},\, \alpha \downarrow \lambda}\mathcal{V}(\alpha),
\end{equation*}
where $\hat{K}$ is the set of dominant Spin($n+1$) weights and $\mathcal{V}(\alpha)$ is the $\alpha$-isotypic component satisfying the classical branching rule of $K$ and $M$:
\begin{equation*}
\alpha\downarrow\lambda\text{ iff }\alpha_1-\lambda_1\in \mathbb{Z}\text{ and }
\begin{cases} \alpha_1\ge\lambda_1\ge\alpha_2\ge \cdots \ge\lambda_l\ge |\alpha_{l+1}|, & \text{n odd}\\
\alpha_1\ge\lambda_1\ge\alpha_2\ge \cdots \ge\lambda_{l-1}\ge \alpha_l\ge |\lambda_l|, &\text{n even}.
\end{cases}
\end{equation*}
The conformal action of $G$ and its infinitesimal representation correspond to those in (\ref{int}).

Let $A=A_{2r}$ be an intertwinor of order $2r$ of the $(\gg,K)$ representation. 
That is, a $K$-map satisfying
the intertwining relation
\begin{equation}\label{rel-1}
A\left(\tilde{\cL}_X+\left(\frac{n}2-r\right)\omega\right)=
\left(\tilde{\cL}_X+\left(\frac{n}2+r\right)\omega\right)A\quad \text{for all }X\in \gg\, \text{ with its conformal factor }\omega,
\end{equation}
where $\tilde{\cL}_X$ is the {\em reduced Lie derivative}, $\tilde{\cL}_X=\cL_X+(l-m)\omega$ on tensors of $\twostac{l}{m}$-type. 

$K$ acts as isometries so every conformal vector field in $\mathfrak{k}$ has zero conformal factor. From now on, we fix the conformal vector field $Y=\sin\r\, \ptl_\r\in \mathfrak{s}$ and  its nonzero  conformal factor $\omega=\cos\r$.
The following lemma compares the Lie derivative and covariant derivative on the spinor-$k$-form bundle $\Sigma\wedge^k$ on the sphere.
\begin{lem}\label{p}
For $\Phi\in\Sigma\wedge^k$, 
\begin{align*}
((\cL_Y-\N_Y)\Phi)_{a_1\cdots a_k}&=k\,\omega\Phi_{a_1\cdots a_k}\, .
\end{align*}
\end{lem}
\begin{proof} Note that for a 1-form $\eta$ and a vector field $X$,
$$
\langle(\cL_Y-\nd_Y)\eta,X\rangle=-\langle\eta,(\cL_Y-\nd_Y)X\rangle,
$$
since $\cL_Y-\nd_Y$ kills scalar functions.
But by the symmetry of the Riemannian connection,
$$
[Y,X]-\nd_YX=-\nd_XY.
$$
We conclude that
$$
(\cL_Y-\nd_Y)\eta=\langle\eta,\nd Y\rangle,
$$
where in the last expression, $\langle\cdot,\cdot\rangle$ is the pairing
of a 1-form with the contravariant part of a $\twostac{1}{1}$-tensor:
$$
((\cL_Y-\nd_Y)\eta)_\l=\eta_\m\nd_\l Y^\m.
$$
Since $Y$ is a conformal vector field, 
$$
(\N Y_\flat)_{\l\m}=(\N Y_\flat)_{(\l\m)}+(\N Y_\flat)_{[\l\m]}=(\omega g+\frac{1}{2}
dY_\flat)_{\l\m}=\omega g_{\l\m}\, .
$$
Thus
\begin{equation*}
((\cL_Y-\N_Y)\eta)_\l=\omega\eta_\l\, .
\end{equation*}
Since $\cL_Y-\N_Y$ is a derivation, for $\Psi\in \wedge^k$,
\begin{equation*}
((\cL_Y-\N_Y)\Psi)_{a_1\cdots a_k}
=\Psi_{\l a_2\cdots a_k}\N_{a_1}Y^\l+\cdots+\Psi_{a_1\cdots a_{k-1}\l}\N_{a_k}Y^\l
=k\,\omega\Psi_{a_1\cdots a_k}\, .
\end{equation*}
On the spinor bundle $\Sigma$, on the other hand, \cite[eq(16)]{Kosmann:72}
$$
\cL_Y-\nd_Y=-\tfrac14\nd_{[a}Y_{b]}\g^a\g^b
=-\tfrac18(dY_\flat)_{ab}\g^a\g^b=0,
$$
where $\gamma$ is a fundamental tensor-spinor. That is, a smooth section of the bundle $TS^n\otimes\text{End}(\Sigma)$ satisfying
$$
\gamma^a\gamma^b+\gamma^b\gamma^a=-2g^{ab}\cdot\text{Id}_\Sigma \quad \text{and} \quad \nabla\gamma=0.
$$
Thus the lemma follows.
\end{proof}
The spectrum generating relation that converts (\ref{rel-1}) is given in the following lemma.
\begin{lem}
On tensor-spinors of any type,
\begin{equation*}
[\N^*\N,\omega]=2\left(\nabla_Y+\dfrac{n}{2}\omega\right)\, ,
\end{equation*}
where [,] is the operator commutator.
\end{lem}
\begin{proof}
If $\varphi$ is any smooth section, then
\begin{equation*}
[\N^*\N,\omega]\varphi=(\Lap \omega)\varphi-2\iota (d\omega) \nabla\varphi
=(n\omega+2\iota(Y)\N)\varphi
=(n\omega+2\N_Y)\, \varphi\, ,
\end{equation*}
where $\iota$ is the interior multiplication.
\end{proof}
Thus the intertwining relation (\ref{rel-1}) on spinor-$k$-forms becomes
\begin{equation}\label{rel-2}
A\left(\frac{1}{2}[\N^*\N,\omega]-r\omega\right)=
\left(\frac{1}{2}[\N^*\N,\omega]+r\omega\right)A\, .
\end{equation}
Now we assume $n \ge 3$ is odd and look at the intertwinors on 
$\mathbb{V}(\underbrace{\tfrac{3}{2},\cdots,\tfrac{3}{2}}_{k},\frac{1}{2}, \cdots,\frac{1}{2})$.
\noindent {\bf Case I}: $k=0$.\\
Let $\cV_\e(j):=\cV(\frac{1}{2}+j,\frac{1}{2},\cdots,\frac{1}{2},\frac{\e}{2})$, where $\e=\pm 1$. Note that a proper conformal factor takes a section of $\cV_\e(j)$ to a direct sum of sections of $\cV_\e(j+1)$, $\cV_\e(j-1)$, and $\cV_{-\e}(j)$ : 
\begin{equation}\label{arrows}
\begin{array}{ccccc}
\cV_\e(j+1) & & & &  \\
&\nwarrow& & &\\
& & \cV_\e(j) &\rightarrow &\cV_{-\e}(j)\\
&\swarrow& & & \\
\cV_\e(j-1)& & & &
\end{array}
\end{equation}
Apply the intertwining relation (\ref{rel-2}) to a section $\varphi$ in $\cV_\e(j)$ :
\begin{equation*}
\begin{array}{c}
A\left(\frac{1}{2}[\N^*\N,\omega]-r\omega\right)\varphi=
\left(\frac{1}{2}[\N^*\N,\omega]+r\omega\right)A\varphi \\
\Leftrightarrow
A\left(\frac{1}{2}(\N^*\N(\omega\varphi)-\omega(\N^*\N_\a\varphi))-r\omega\varphi\right)=\mu_\a
\left(\frac{1}{2}(\N^*\N(\omega \varphi)-\omega(\N^*\N_\a\varphi))+r\omega \varphi\right)\, ,
\end{array}
\end{equation*}
where $\mu_\a$ (resp. $\N^*\N_\a$) is the eigenvalue of $A$ (resp. $\N^*\N$) on the $K$-type $\a:=\cV_\e(j)$.
Let Proj${}_\beta\omega|_\a\varphi$ be the projection of $\omega\varphi$ onto the $K$-type $\b$ direct summand. 
The ``compression'', from the $K$-type $\a$ to the $K$-type $\beta$, 
of the above relation becomes \underline{Proj${}_\beta\omega|_\a\varphi$ times}
\begin{equation}\label{numer}
\left(\frac12\N^*\N|^\beta_\a+r\right)\mu_\a
=\left(\frac12\N^*\N|^\beta_\a-r\right)\mu_\b\,,  
\end{equation}
where $\mu_\beta$ (resp. $\N^*\N_\beta$) is the eigenvalue of $A$ (resp. $\N^*\N$) on
the $K$-type $\beta$ and $\N^*\N|^\beta_\alpha:=\N^*\N_\beta-\N^*\N_\alpha$.  The underlined phrase above is a key point. We
have achieved a factorization in which one factor is purely numerical
(that appearing in (\ref{numer})).  ``Canceling'' the other factor,
Proj${}_\beta\omega|_\alpha$, we get purely numerical recursions that 
are guaranteed to give intertwinors.  If we wish to see the {\em uniqueness}
of intertwinors this way, we need to establish the nontriviality of 
the Proj${}_\beta\omega|_\alpha$. In fact this nontriviality follows from
Branson \cite[sect. 6]{Branson:97-1}.

Let $J=\frac{n}{2}+j$. The transition quantities $\m_\b/\m_\a$ with respect to the diagram (\ref{arrows}) are
\begin{equation*}
\begin{array}{ccccc}
(J+\frac{1}{2}+r)/(J+\frac{1}{2}-r) & & & &  \\
&\nwarrow& & &\\
& &\bullet &\rightarrow &-1\\
&\swarrow& & & \\
(-J+\frac{1}{2}+r)/(-J+\frac{1}{2}-r)& & & &
\end{array}
\end{equation*}
Choosing normalization $\m_{\cV_1(0)}=1$, we get 
\begin{thm}
The unique spectral function $Z_\e(r;j)$ on $\cV_\e(j)$ is, up to normalization, 
\begin{equation*}
Z_\e(r;j)=\e\cdot \frac{\Gamma(J+\frac{1}{2}+r)\Gamma(\frac{n}{2}+\frac{1}{2}-r)}{\Gamma(J+\frac{1}{2}-r)\Gamma(\frac{n}{2}+\frac{1}{2}+r)}\, .
\end{equation*}
\end{thm}
Note that when $2r=1$, $Z_\e(1/2;j)$ is a constant multiple of the Dirac operator $\dc=\gamma^a\nabla_a$. In general, we have 
\begin{cor}\cite{BO:06}. \label{odd}
For l=0,1,2,$\cdots$, 
\begin{equation*}
D_{2l+1}:=\dc(\dc^2-1^2)\cdots(\dc^2-l^2)
\end{equation*}
is a differential intertwinor of order $2r=2l+1$.
\end{cor}
\begin{proof}
Clearly, $Z_\e((2l+1)/2;j)$ and $D_{2l+1}$ differ by a nonzero constant.
\end{proof}
\noindent {\bf Case II}: $1\le k \le (n-1)/2$.\\
Let, for $q=0,1$ and $\e=\pm 1$, 
\begin{equation*}
\cV_\e(j,q):=
\left\{\begin{array}{l}
\cV(\frac{3}{2}+j,\underbrace{\tfrac{3}{2},\cdots,\tfrac{3}{2}}_{k-1},\frac{1}{2}+q,\frac{1}{2},\cdots,\frac{1}{2},\frac{\e}{2}),\quad 1\le k<(n-1)/2,\\
\cV(\frac{3}{2}+j,\tfrac{3}{2},\cdots,\tfrac{3}{2},\e\cdot(\frac{1}{2}+q),\quad k=(n-1)/2.
\end{array}\right.
\end{equation*}

We have a diagram of direct summand $K$-types centered at $\cV_\e(j,0)$ :
\begin{equation*}
\begin{array}{ccccc}
 & & \cV_\e(j+1,0) & &  \\
& & \uparrow & &\\
\cV_\e(j,1) & \leftarrow & \cV_\e(j,0) &\rightarrow &\cV_{-\e}(j,0)\\
& & \downarrow & & \\
& & \cV_\e(j-1,0) & &
\end{array}
\end{equation*}

Let $L=\frac{n}{2}+1+j$. The transition quantities with respect to the above diagram are
\begin{equation*}
\begin{array}{ccc}
  & (L+\frac{1}{2}+r)/(L+\frac{1}{2}-r) &  \\
& \uparrow &\\
(\frac{n}{2}-k+\frac{1}{2}+r)/(\frac{n}{2}-k+\frac{1}{2}-r) & \leftarrow \quad  \bullet  \quad \rightarrow &-1\\
& \downarrow & \\
& (-L+\frac{1}{2}+r)/(-L+\frac{1}{2}-r) &
\end{array}
\end{equation*}
With respect to the following diagram 
\begin{equation*}
\begin{array}{ccccc}
 & & \cV_\e(j+1,1) & &  \\
& & \uparrow & &\\
\cV_\e(j,0) & \leftarrow & \cV_\e(j,1) &\rightarrow &\cV_{-\e}(j,1)\\
& & \downarrow & & \\
& & \cV_\e(j-1,1) & &
\end{array}
\end{equation*}
we get the transition quantities
\begin{equation*}
\begin{array}{ccc}
  & (L+\frac{1}{2}+r)/(L+\frac{1}{2}-r) &  \\
& \uparrow &\\
(-\frac{n}{2}+k-\frac{1}{2}+r)/(-\frac{n}{2}+k-\frac{1}{2}-r) & \leftarrow \quad  \bullet \quad \rightarrow &-1\\
& \downarrow & \\
& (-L+\frac{1}{2}+r)/(-L+\frac{1}{2}-r) &
\end{array}
\end{equation*}
Choosing normalization $\m_{\cV_1(0,1)}=1=(\frac{n}{2}-k+\frac{1}{2}-r)/(\frac{n}{2}-k+\frac{1}{2}+r)\cdot\m_{\cV_1(0,0)}$, we get
\begin{thm}
The unique spectral functions $Z_\e(r,j,1)$ and $Z_\e(r,j,0)$ on $\cV_\e(j,1)$ and $\cV_\e(j,1)$ respectively are, up to normalization,
\begin{align*}
&Z_\e(r,j,1)=\e\cdot \frac{\Gamma(L+\frac{1}{2}+r)\Gamma(\frac{n}{2}+\frac{3}{2}-r)}{\Gamma(L+\frac{1}{2}-r)\Gamma(\frac{n}{2}+\frac{3}{2}+r)} \quad \text{and}\\
&Z_\e(r,j,0)=\e\cdot\frac{n-2k+1-2r}{n-2k+1+2r}\cdot \frac{\Gamma(L+\frac{1}{2}+r)\Gamma(\frac{n}{2}+\frac{3}{2}-r)}{\Gamma(L+\frac{1}{2}-r)\Gamma(\frac{n}{2}+\frac{3}{2}+r)}\, .
\end{align*}
\end{thm}
When $2r=1$ and $k=1$, up to a constant, $Z_\e(1/2,j,1)=\e\cdot L$ and $Z_\e(1/2,j,0)=\e\cdot \frac{n-2}{n}\cdot L$. The Rarita-Schwinger operator 
\begin{equation*}
(S\varphi)_a:=\gamma^b\N_b\varphi_a-\frac{2}{n}\gamma_a\N^b\varphi_b
\end{equation*}
eigenvalues agree with the above data.

For $\varphi\in\Sigma\wedge^k$, we define the following convenient operators \cite{Branson:99}:
\begin{align*}
(\tilde{d}\varphi)_{a_0\cdots a_k} &:=\sum_{i=0}^k (-1)^i\N_{a_i}\varphi_{a_0\cdots a_{i-1}a_{i+1}\cdots a_k},\\
(\tilde{\delta}\varphi)_{a_2\cdots a_k}&:=-\N^b\varphi_{ba_2\cdots a_k},\\
(\e(\gamma)\varphi)_{a_0\cdots a_k}&:=\sum_{i=0}^k (-1)^i\gamma_{a_i}\varphi_{a_0\cdots a_{i-1}a_{i+1}\cdots a_k}, \\
(\iota(\gamma)\varphi)_{a_2\cdots a_k}&:=\gamma^b\varphi_{ba_2\cdots a_k},\\
(\mathbb{D}\varphi)_{a_1\cdots a_k}&:=(\iota(\gamma)\tilde{d}+\tilde{d}\iota(\gamma))\varphi)_{a_1\cdots a_k}=-(\tilde{\delta}\e(\gamma)+\e(\gamma)\tilde{\delta})\varphi)_{a_1\cdots a_k}=\gamma^\l\N_\l\varphi_{a_1\cdots a_k}.
\end{align*}
The operator 
\begin{equation*}
P_k:=\frac{n-2k+4}{2}\iota(\gamma)\tilde{d}+\frac{n-2k}{2}(\tilde{d}\iota(\gamma)-\tilde{\delta}\e(\gamma))-\frac{n-2k-4}{2}\e(\gamma)\tilde{\delta}
\end{equation*}
is conformally covariant on $\Sigma\wedge^k$ and its restriction to $\mathbb{T}^k:=\mathbb{V}(\underbrace{\tfrac{3}{2},\dots\tfrac{3}{2}}_{k},\frac{1}{2},\dots\frac{1}{2})$
\begin{equation*}
P_k|_{\mathbb{T}^k}=(n-2k+2)\mathbb{D}+2\e(\gamma)\tilde{\delta}
\end{equation*}
is an intertwinor on $\mathbb{T}^k$ with eigenvalues \cite[Example 10.6]{Branson:99}
\begin{equation*}
\e\cdot 2 \cdot L\cdot(q+\frac{n}{2}-k)\, \text{ on }\cV_\e(j,q)\, \text{ for }q=0,1\, .
\end{equation*}
Let us take a convenient normalization of the operator: 
\begin{equation*}
A_{k,0}:=\frac{1}{n-2k+2}\cdot P_k|_{\mathbb{T}^k}
=\left\{\begin{array}{cc} \e\cdot L & \text{on }\cV_\e(j,1),\\ \e\cdot {\displaystyle \frac{n-2k}{n-2k+2}}\cdot L & \text{on }\cV_\e(j,0)\, . \end{array}\right.
\end{equation*}
Note that $A_{0,0}=\dc$, the Dirac operator and $A_{1,0}=S$, the Rarita-Schwinger operator.
 
Consider the operator $T_{k-1}: \mathbb{T}^{k-1}\rightarrow\mathbb{T}^k$ defined by 
\begin{equation*}
T_{k-1}=\frac{1}{k}\tilde{d}+\frac{1}{k(n-2(k-1))}\e(\gamma)\mathbb{D}+\frac{1}{k(n-2(k-1))(n-2(k-1)+1)}\e(\gamma)^2\tilde{\delta}\, .
\end{equation*}
This is the orthogonal projection of $\nabla$ onto $\mathbb{T}^k$ summand ($1/k\cdot\tilde{d}^{\text{top}}_{k-1}$ in \cite{BH:02}):
\begin{equation}\label{tk}
\mathbb{T}^{k-1}\stackrel{\nabla}{\longrightarrow}T^*S^n\otimes\mathbb{T}^{k-1}\cong_{\text{Spin}(n)}
\left\{\begin{array}{l}
\mathbb{T}^{k-2}\oplus\mathbb{T}^{k-1}\oplus\mathbb{T}^k\oplus\mathbb{Z}^{k-1},\quad 1\le k <(n-1)/2,\\
\mathbb{T}^{k-2}\oplus\mathbb{T}^{k-1}\oplus\mathbb{T}^k, \quad k=(n-1)/2,
\end{array}\right.
\end{equation}
where $\mathbb{Z}^{k-1}\cong_{\text{Spin($n$)}}\mathbb{V}(\tfrac{5}{2},\underbrace{\tfrac{3}{2},\dots,\tfrac{3}{2}}_{k-2},\tfrac{1}{2},\dots\tfrac{1}{2})$.
Note also that the formal adjoint of $T_{k-1}$ is $T^*_{k-1}=\tilde{\delta}$. 
\begin{lem}
The second order operator $T_{k-1}T_{k-1}^*$ acts as a scalar
\begin{equation*}
\begin{array}{cl}
0 &\quad \text{ on }\cV_\epsilon (j,1)\, \text{ and}\\ 
\displaystyle{\frac{(n-2k+1)(L^2-(n/2-k+1)^2)}{k(n-2k+2)}} &\quad \text{ on }\cV_\epsilon (j,0)\, . 
\end{array}
\end{equation*}
\end{lem}
\begin{proof}
Let $G_{\mathbb{T}^{k-1}}$ be the gradient from $\mathbb{T}^{k}$ to a copy of $\mathbb{T}^{k-1}$. That is, the orthogonal projection of $\nabla$ onto $\mathbb{T}^{k-1}$ summand from $\mathbb{T}^{k}$ (See (\ref{tk})). And let $G_{\mathbb{T}^{k-1}}^*$  be the formal adjoint of $G_{\mathbb{T}^{k-1}}$. By theorem 8.3 in \cite{Branson:99}, $G_{\mathbb{T}^{k-1}}^*G_{\mathbb{T}^{k-1}}$ over $\mathbb{T}^k$ acts as a scalar 
\begin{equation*}
\begin{array}{cl}
0&\quad\text{ on }\cV_\epsilon (j,1),\\
\displaystyle{\frac{(n-2k+3)(L^2-(n/2-k+1)^2)}{(n-k+2)(n-2k+2)}} &\quad\text{ on }\cV_\epsilon (j,0)
\end{array}
\end{equation*}
and $T_{k-1}^*T_{k-1}$ over $\mathbb{T}^{k-1}$ acts as a scalar
\begin{equation*}
\begin{array}{cl}
\lambda:=\displaystyle{\frac{(n-2k+1)(L^2-(n/2-k+1)^2)}{k(n-2k+2)}} &\quad\text{ on }\cV_\epsilon (j,1),\\
0&\quad\text{ on }\cV_\epsilon (j,0)\, .
\end{array}
\end{equation*}
Note that $T_{k-1}^*$ is a first order differential operator from $\mathbb{T}^{k}$ to $\mathbb{T}^{k-1}$. By uniqueness of the gradient up to a constant \cite{F:76}, $T_{k-1}T_{k-1}^*$ is a scalar multiple of $G_{\mathbb{T}^{k-1}}^*G_{\mathbb{T}^{k-1}}$. In particular, $T_{k-1}T_{k-1}^*=0$ on $\cV_\epsilon (j,1)$. To determine $T_{k-1}T_{k-1}^*$ on $\cV_\epsilon (j,0)$, we take $\varphi\in \cV_\epsilon(j,1)$ over $\mathbb{T}^{k-1}$. Then, $T_{k-1}\varphi\in \cV_\epsilon(j,0)$ over $\mathbb{T}^k$. Thus, $T_{k-1}T_{k-1}^*(T_{k-1}\varphi)=\lambda\cdot T_{k-1}\varphi$.
\end{proof}
\begin{rmk}
The eigenvalue formulas in the lemma are valid for both odd and even dimensional cases. They also imply that $\cV(j,0)$ summands make up $\mathcal{R}(T_{k-1})$ and $\cV(j,1)$ summands make up $\mathcal{R}(T_k^*)$ over $\mathbb{T}^k$ for $1\le k \le (n-2)/2$, where $\mathcal{R}(T_{k-1})$ (resp. $\mathcal{R}(T_k^*)$) denotes the range of $T_{k-1}$ (resp. $T_k^*$).
\end{rmk}
Thus, for $l=1,2,\cdots$, the operator
\begin{equation*}
A_{k,l}:=A_{k,0}^2-l^2\cdot\text{Id}+a_l\cdot T_{k-1}T_{k-1}^*,\\
\end{equation*}
where 
\begin{equation*}
a_l=\frac{-16kl^2}{(n-2k+2)(n-2k+2-2l)(n-2k+2+2l)}
\end{equation*}
acts as a scalar
\begin{equation*}
\begin{array}{cl}
L^2-l^2 & \text{ on }\cV_\e(j,1)\, \text{ and}\\
{\displaystyle \frac{(n-2k-2l)(n-2k+2l)}{(n-2k+2-2l)(n-2k+2+2l)}}\cdot(L^2-l^2)&\text{ on }\cV_\e(j,0).
\end{array}
\end{equation*}
Collecting the above information, we get 
\begin{thm}
The operator $\displaystyle{\prod_{a=0}^l A_{k,a}}$ for $l=1,2,\cdots$, acts as a scalar 
\begin{align*}
&\e\cdot L(L^2-1^2)\cdots (L^2-l^2) \quad \text{on }\cV_\e(j,1)\, \text{ and}\\
&\e\cdot \frac{n-2k-2l}{n-2k+2+2l}\cdot L(L^2-1^2)\cdots(L^2-l^2) \quad \text{on }\cV_\e(j,0).
\end{align*}
Consequently, it is a differential intertwinor of order 2l+1 on $\mathbb{T}^k$.
\end{thm}
\begin{proof}
$Z_\e((2l+1)/2,j,1)$ is a nonzero constant multiple of  $\displaystyle{\prod_{a=0}^l A_{k,a}}$ on $\cV_\e(j,1)$ and $Z_\e((2l+1)/2,j,0)= \dfrac{n-2k-2l}{n-2k+2+2l}\cdot Z_\e((2l+1)/2,j,1)$. 
\end{proof}

\vspace{1cm}
\noindent
Department of Mathematics\\
University of North Dakota\\
Grand Forks ND 58202 USA\\
email: doojin.hong@email.und.edu

\end{document}